\documentclass[12pt]{article}
\usepackage{graphics}
\usepackage{float}
\usepackage{amsfonts}
\usepackage{amsmath}

\input{tcilatex}
\begin{document}

\author{Ozlem Ersoy and Idris Dag \\
Eski\c{s}ehir Osmangazi \"{U}niversity, Faculty of Science and Art,\\
Department of Mathematics-Computer, Eski\c{s}ehir, Turkey,}
\title{\textbf{The Numerical Approach to the Fisher's Equation via
Trigonometric Cubic B-spline Collocation Method}}
\maketitle

\begin{abstract}
In this study, we set up a numerical technique to get approximate solutions
of Fisher's equation which is one of the most important model equation in
population biology. We integrate the equation fully by using combination of
the trigonometric cubic B-spline functions for space variable and
Crank-Nicolson for the time integration. Numerical results have been
presented to show the accuracy of the current algorithm. We have seen that
the proposed technique is a good alternative to some existing techniques for
getting solutions of the Fisher's equation.

\end{abstract}

\section{Introduction}

One of the most interesting equation in physical phenomena is
reaction-diffusion equation. We focus on one of the special case of the
reaction-diffusion given by%
\begin{equation}
U_{t}=\lambda U_{xx}+\beta U(1-U),\text{ }-\infty <x<\infty  \label{f2}
\end{equation}%
where $\beta $ is a real parameter. This equation is known as a Fisher's
equation(FE) which was introduced by Fisher \cite{first} to describe the
kinetic advancing rate of an advantageous gene.\ In a large number of
biological and chemical phenomena, the reaction term is represented by $%
\beta U(1-U)$, where $\beta >0$ can be dependent of the space variable. FE
represents the evolution of the population due to the two competing physical
processes and changes of interaction of diffusion and nonlinear reaction can
be observed.

The trigonometric cubic B-spline (TCB) functions have been started to adapt
the numerical methods for obtaing the solutions of the differential
equations although use of these are not widespread as the approximate
functions in the numerical methods. The numerical methods for solving a type
of ordinary differential equations with trigonometric quadratic and cubic
splines are given by A. Nikolis in the papers\cite{g1,g2}. The linear
two-point boundary value problems of order two are solved using TCB
interpolation method\cite{we}. The another numerical method employed the TCB
is set up to solve a class of linear two-point singular boundary value
problems in the study\cite{w1}. Very recently a collocation finite
difference scheme based on TCB is developed for the numerical solution of a
one-dimensional hyperbolic equation (wave equation) with non-local
conservation condition\cite{aa}. A new two-time level implicit technique
based on the TCB is proposed for the approximate solution of the
nonclassical diffusion problem with nonlocal boundary condition in the study%
\cite{ma}.

The numerical approaches for types of the Fisher equation have been given to
model some events, which can not be exhibited by theoretical means. Although
various numerical tecniques are employed for getting numerical solutions of
the Fisher equation, we only documents the spline-based numerical methods.
Orthogonal cubic spline collocation method is constructed for the extended
Fisher-Kolmogorov equation\cite{fs9}. The Galerkin method constructed using
quartic B-splines and the Crack-Nicolson method is employed for the
numerical solution of FE\cite{fs4}. A numerical solution of the extended
Fisher-Kolmogorov equation by using the quintic B-spline collocation scheme
is proposed\cite{fs8}. FE is fully discretized by the finite element method
based on the quadratic B-spline Galerkin technique in space and by use of
the Crank-Nicolson method for the integration of the obtained matrix
ordinary differental equations in space\cite{fs5}. Both the cubic B-spline
collocation method and the usual finite difference tecniques is applied for
getting the solutions of the well-known Fisher's equation which was
discretized in space and time respectively\cite{fs3}. The Cubic B-spline
quasi-interpolation is presented, by using the derivative of the
quasi-interpolation to approximate the spatial derivative of the dependent
variable and a low order forward difference to approximate the time
derivative of the dependent variable for finding solutions of the
Burgers--Fisher equation\cite{fs6}. Combination of the quintic B-spline
based collocation method and crank-Nicolson tecnique is designed to
discretize the Fisher's equation fully \cite{fs7}. Crank-Nicolson method is
used for time discretization. \ A collocation method based on modified cubic
B-splines for spatial variable and derivatives is applied to the Fisher's
equation to produce a system of first-order ordinary differential equations,
Which was solved by SSP-RK54 scheme\cite{fs2}. The spline approximation has
been used in spatial direction together with the finite difference
approximation in time \ have been set up to gave numerical solution of the
Fisher equation\cite{fs1}. The exponential cubic B-spline collocation method
is constructed to obtain numerical solutions of the Fisher equation\cite%
{fs10}.

The purpose of this studying is to obtain the numerical solutions of the
Fisher's equation via the trigonometric cubic B-spline collocation algorithm
and research the solvability of the reaction-diffusion type equations by the
B-spline finite element methods. Crank-Nicolson formulas and trigonometric
cubic B-spline functions are used for time and space discretization
respectively. Over the uniform mesh, Crank-Nicolson formulas are employed
for time discretization whereas Rubin and Graves\cite{rubin} technique is
used for the linearization. The collocation method together with
trigonometric cubic B-spline approximations exhibits an economical
alternative since it only requires the evaluation of the unknown parameters
at the grid points.

The different initial and boundary conditions (BCs) for Fisher's Equation(%
\ref{f2}) are given as

\begin{eqnarray}
u(x,0) &=&u_{x}(0)\in \lbrack 0.1],\text{ }x\in \lbrack -\infty ,\infty ]
\label{s1} \\
\lim_{x\rightarrow -\infty }u(x,t) &=&1,\lim_{x\rightarrow \infty }u(x,t)=0
\label{s2} \\
\lim_{x\rightarrow \pm \infty }u(x,t) &=&0  \label{s3}
\end{eqnarray}

In the literature, conditions (\ref{s1}) and (\ref{s2}) together are
commonly referred as nonlocal conditions,while conditions (\ref{s1}) and (%
\ref{s3}) together are usually referred as local conditions.

\section{Cubic Trigonometric B-spline Collocation Method}

Consider a uniform partition of the problem domain $[a=x_{0},b=x_{N}]$ at
the knots $x_{i},i=0,...,N$ with mesh spacing $h=(b-a)/N.$ On this partition
together with additional knots $x_{-1},x_{0},x_{N+1},x_{N+2},x_{N+3}$
outside the problem domain, $CTB_{i}${}$(x)$ can be defined as

\begin{equation}
CTB_{i}(x)=\frac{1}{\theta }\left \{ 
\begin{tabular}{ll}
$\omega ^{3}(x_{i-2}),$ & $x\in \left[ x_{i-2},x_{i-1}\right] $ \\ 
$\omega (x_{i-2})(\omega (x_{i-2})\phi (x_{i})+\phi (x_{i+1})\omega
(x_{i-1}))+\phi (x_{i+2})\omega ^{2}(x_{i-1}),$ & $x\in \left[ x_{i-1},x_{i}%
\right] $ \\ 
$\omega (x_{i-2})\phi ^{2}(x_{i+1})+\phi (x_{i+2})(\omega (x_{i-1})\phi
(x_{i+1})+\phi (x_{i+2})\omega (x_{i})),$ & $x\in \left[ x_{i},x_{i+1}\right]
$ \\ 
$\phi ^{3}(x_{i+2}),$ & $x\in \left[ x_{i+1},x_{i+2}\right] $ \\ 
$0,$ & $\text{otherwise}$%
\end{tabular}%
\right.  \label{r2}
\end{equation}%
where $\omega (x_{i})=\sin (\frac{x-x_{i}}{2}),\phi (x_{i})=\sin (\frac{%
x_{i}-x}{2}),\theta =\sin (\frac{h}{2})\sin (h)\sin (\frac{3h}{2}).$

$CTB_{i}(x)$ are twice continuously differentiable piecewise trigonometric
cubic B-spline on the interval$[a,b]$. The iterative formula \bigskip

\begin{equation}
T_{i}^{k}(x)=\frac{\sin (\frac{x-x_{i}}{2})}{\sin (\frac{x_{i+k-1}-x_{i}}{2})%
}T_{i}^{k-1}(x)+\frac{\sin (\frac{x_{i+k}-x}{2})}{\sin (\frac{x_{i+k}-x_{i+1}%
}{2})}T_{i+1}^{k-1}(x),\text{ }k=2,3,4,...  \label{r3}
\end{equation}%
gives the cubic B-spline trigonometric functions starting with the
CTB-splines of order $1$

\begin{equation*}
T_{i}^{1}(x)=\left \{ 
\begin{tabular}{l}
$1,\ x\in \lbrack x_{i},x_{i+1})$ \\ 
$0$, otherwise.%
\end{tabular}%
\right.
\end{equation*}

Each $CTB_{i}(x)$ is twice continuously differentiable and the values of $%
CTB_{i}(x),CTB_{i}^{^{\prime }}(x)$ and $CTB_{i}^{^{\prime \prime }}(x)$ at
the knots $x_{i}$ 's can be computed from Eq.(\ref{r3}) as

\begin{equation*}
\begin{tabular}{l}
Table 1: Values of $B_{i}(x)$ and its principle two \\ 
derivatives at the knot points \\ 
\begin{tabular}{|l|l|l|l|}
\hline
& $T_{i}(x_{k})$ & $T_{i}^{\prime }(x_{k})$ & $T_{i}^{\prime \prime }(x_{k})$
\\ \hline
$x_{i-2}$ & $0$ & $0$ & $0$ \\ \hline
$x_{i-1}$ & $\sin ^{2}(\frac{h}{2})\csc \left( h\right) \csc (\frac{3h}{2})$
& $\frac{3}{4}\csc (\frac{3h}{2})$ & $\frac{3(1+3\cos (h))\csc ^{2}(\frac{h}{%
2})}{16\left[ 2\cos (\frac{h}{2})+\cos (\frac{3h}{2})\right] }$ \\ \hline
$x_{i}$ & $\frac{2}{1+2\cos (h)}$ & $0$ & $\frac{-3\cot ^{2}(\frac{3h}{2})}{%
2+4\cos (h)}$ \\ \hline
$x_{i+1}$ & $\sin ^{2}(\frac{h}{2})\csc \left( h\right) \csc (\frac{3h}{2})$
& $-\frac{3}{4}\csc (\frac{3h}{2})$ & $\frac{3(1+3\cos (h))\csc ^{2}(\frac{h%
}{2})}{16\left[ 2\cos (\frac{h}{2})+\cos (\frac{3h}{2})\right] }$ \\ \hline
$x_{i+2}$ & $0$ & $0$ & $0$ \\ \hline
\end{tabular}%
\end{tabular}%
\end{equation*}%
$CTB_{i}(x)$ , $i=-1,...,N+1$ are a basis for the trigonometric spline
space. An approximate solution $U_{N}$ to the unknown $U$ is written in
terms of the expansion of the CTB as

\begin{equation}
U_{N}(x,t)=\sum_{i=-1}^{N+1}\delta _{i}CTB_{i}(x)  \label{r4}
\end{equation}%
where $\delta _{i}$ are time dependent parameters to be determined from the
collocation points $x_{i},$ $i=0,...,N$ and the boundary and initial
conditions. The nodal values $U$ and its first and second derivatives at the
knots can be found from the (\ref{r4}) as 
\begin{equation}
\begin{tabular}{l}
$U_{i}=\alpha _{1}\delta _{i-1}+\alpha _{2}\delta _{i}+\alpha _{1}\delta
_{i+1}$ \\ 
$U_{i}^{\prime }=\beta _{1}\delta _{i-1}+\beta _{2}\delta _{i+1}$ \\ 
$U_{i}^{\prime \prime }=\gamma _{1}\delta _{i-1}+\gamma _{2}\delta
_{i}+\gamma _{1}\delta _{i+1}$%
\end{tabular}
\label{r5}
\end{equation}%
\begin{equation*}
\begin{array}{ll}
\alpha _{1}=\sin ^{2}(\frac{h}{2})\csc (h)\csc (\frac{3h}{2}) & \alpha _{2}=%
\dfrac{2}{1+2\cos (h)} \\ 
\beta _{1}=-\frac{3}{4}\csc (\frac{3h}{2}) & \beta _{2}=\frac{3}{4}\csc (%
\frac{3h}{2}) \\ 
\gamma _{1}=\dfrac{3((1+3\cos (h))\csc ^{2}(\frac{h}{2}))}{16(2\cos (\frac{h%
}{2})+\cos (\frac{3h}{2}))} & \gamma _{2}=-\dfrac{3\cot ^{2}(\frac{h}{2})}{%
2+4\cos (h)}%
\end{array}%
\end{equation*}

The Crank--Nicolson \ scheme is used to discretize time variables of the
unknown $U$ in the Fisher's equation so that one obtain the time discretized
form of the equation as%
\begin{equation}
\frac{U^{n+1}-U^{n}}{\Delta t}=\lambda \frac{(U_{xx})^{n+1}+(U_{xx})^{n}}{2}%
+\beta \frac{U^{n+1}+U^{n}}{2}-\beta \frac{(UU)^{n+1}+(UU)^{n}}{2}
\label{d1}
\end{equation}%
where $U^{n+1}=U(x,t)$ is the solution of the equation at the $(n+1)$th time
level. Here $t^{n+1}$ $=t^{n}+t$, and $\Delta t$ is the time step,
superscripts denote $n$ th time level , $t^{n}=n\Delta t.$ The nonlinear
term $(U^{2})^{n+1}$ in Eq. (\ref{d1}) may be linearized by using the
following term \cite{rubin}:%
\begin{equation}
(U^{2})^{n+1}=2U^{n}U^{n+1}-(U^{n})^{2}  \label{d2}
\end{equation}%
we get%
\begin{equation}
\frac{2}{\Delta t}U^{n+1}-\lambda U_{xx}^{n+1}-\beta U^{n+1}+2\beta KU^{n+1}=%
\frac{2}{\Delta t}U^{n}+\lambda U_{xx}^{n}+\beta U^{n}  \label{d3}
\end{equation}

Substitution of (\ref{r4}) into (\ref{d3}) leads to the fully-discretized
equation:

\begin{equation}
\chi _{1}\delta _{m-1}^{n+1}+\chi _{2}\delta _{m}^{n+1}+\chi _{1}\delta
_{m+1}^{n+1}=\chi _{3}\delta _{m-1}^{n}+\chi _{4}\delta _{m}^{n}+\chi
_{3}\delta _{m+1}^{n}  \label{d4}
\end{equation}%
where%
\begin{eqnarray*}
\chi _{1} &=&\left( \frac{2}{\Delta t}-\beta +2\beta K\right) \alpha
_{1}-\lambda \gamma _{1} \\
\chi _{2} &=&\left( \frac{2}{\Delta t}-\beta +2\beta K\right) \alpha
_{2}-\lambda \gamma _{2} \\
\chi _{3} &=&\left( \frac{2}{\Delta t}+\beta \right) \alpha _{1}+\lambda
\gamma _{1} \\
\chi _{4} &=&\left( \frac{2}{\Delta t}+\beta \right) \alpha _{2}+\lambda
\gamma _{2} \\
K &=&\alpha _{1}\delta _{i-1}+\delta _{i}+\alpha _{1}\delta _{i+1}
\end{eqnarray*}%
\begin{equation*}
\begin{array}{ll}
\alpha _{1}=\sin ^{2}(\frac{h}{2})\csc (h)\csc (\frac{3h}{2}) & \alpha _{2}=%
\dfrac{2}{1+2\cos (h)} \\ 
\gamma _{1}=\dfrac{3((1+3\cos (h))\csc ^{2}(\frac{h}{2}))}{16(2\cos (\frac{h%
}{2})+\cos (\frac{3h}{2}))} & \gamma _{2}=-\dfrac{3\cot ^{2}(\frac{h}{2})}{%
2+4\cos (h)}%
\end{array}%
\end{equation*}

The system consist of $N+1$ linear equation in $N+3$ unknown parameters $%
\mathbf{d}^{n+1}=(\delta _{-1}^{n+1},\delta _{0}^{n+1},\ldots ,\delta
_{N+1}^{n+1})$. To make solvable the system, boundary conditions $\sigma
_{1}=U_{0},$ $\sigma _{2}=U_{n}$ are used to find two additional linear
equations: 
\begin{eqnarray}
\delta _{-1} &=&\frac{1}{\alpha _{1}}\left( U_{0}-\alpha _{2}\delta
_{0}-\alpha _{3}\delta _{1}\right) ,  \label{f10} \\
\delta _{N+1} &=&\frac{1}{\alpha _{3}}\left( U_{n}-\alpha _{1}\delta
_{N-1}-\alpha _{2}\delta _{N}\right) .  \notag
\end{eqnarray}%
(\ref{f10}) can be used to eliminate $\delta _{-1},\delta _{N+1}$ from the
system (\ref{d4}) which then becomes the solvable matrix equation for the
unknown $\delta _{0}^{n+1},\ldots ,\delta _{N}^{n+1}.$ A tridiogonal system
of equation can be solved with tribanded Thomas algorithm.

Initial parameters $\delta _{-1}^{0},\delta _{0}^{0},\ldots ,\delta
_{N+1}^{0}$ can be determined from the initial condition and first space
derivative of the initial conditions at the boundaries as the following:

\begin{enumerate}
\item $U_{N}(x_{i},0)$ $=U(x_{i},0),$ $i=0,...,N$

\item $(U_{x})_{N}(x_{0},0)=U^{\prime }(x_{0})(U_{x})_{N}(x_{N},0)=U^{\prime
}(x_{N}).$
\end{enumerate}

\section{Numerical tests}

Numerical method described in the previous section is tested on three
problems for getting solution of \ the Fisher's equation and on one problem
for getting the solution of the Fisher equation in order to demonstrate the
robustness and numerical accuracy.

The discrete error norms $L_{\infty }$%
\begin{equation*}
L_{\infty }=\left \vert U-U_{N}\right \vert _{\infty }=\max \limits_{j}\left
\vert U_{j}-(U_{N})_{j}^{n}\right \vert
\end{equation*}%
is used to measure error between the analytical and numerical solutions. The
relative error%
\begin{equation*}
\text{Relative error}=\sqrt{\frac{\dsum \limits_{j=1}^{N}\left \vert
(U_{N})_{j}^{n+1}-(U_{N})_{j}^{n}\right \vert ^{2}}{\dsum
\limits_{j=1}^{N}\left \vert (U_{N})_{j}^{n}\right \vert ^{2}}}
\end{equation*}%
is used when the analytical solutions does not exit.

(a) Analytical solution of the Fisher equation

\begin{equation*}
u(x,t)=\left( 1+\exp (\sqrt{\frac{\beta }{6}}x-\frac{5\beta }{6}t)\right)
^{-2}
\end{equation*}%
is used in the numerical studies \cite{k1,k2}\ together with the nonlocal
conditions. Over the interval$[a,b]=[-0.2,0.8]$ , the calculation is done
with the number of knots $N=40$ with time step $\Delta t=0.0001$ for
different values of $\beta =2000,$ $5000$ and $10000$ to compare results
with other works \cite{fs3}. The computed results are represented grafically
for some times seen in Figs 1-3.

\begin{equation*}
\begin{array}{c}
\begin{array}{cc}
\begin{tabular}{l}
\FRAME{itbpF}{2.693in}{2.0366in}{0in}{}{}{fig1.jpg}{\special{language
"Scientific Word";type "GRAPHIC";display "USEDEF";valid_file "F";width
2.693in;height 2.0366in;depth 0in;original-width 6.4999in;original-height
4.9165in;cropleft "0";croptop "1";cropright "1";cropbottom "0";filename
'Fig1.jpg';file-properties "XNPEU";}} \\ 
Fig. 1: solutions for $\beta =2000,N=40$%
\end{tabular}
& 
\begin{tabular}{l}
\FRAME{itbpF}{2.6749in}{2.0453in}{0in}{}{}{fig2.jpg}{\special{language
"Scientific Word";type "GRAPHIC";maintain-aspect-ratio TRUE;display
"USEDEF";valid_file "F";width 2.6749in;height 2.0453in;depth
0in;original-width 6.4999in;original-height 4.9562in;cropleft "0";croptop
"1";cropright "1";cropbottom "0";filename 'Fig2.jpg';file-properties
"XNPEU";}} \\ 
Fig. 2: solutions for $\beta =5000,N=40$%
\end{tabular}%
\end{array}
\\ 
\begin{tabular}{l}
\FRAME{itbpF}{2.693in}{2.0781in}{0in}{}{}{fig3.jpg}{\special{language
"Scientific Word";type "GRAPHIC";maintain-aspect-ratio TRUE;display
"USEDEF";valid_file "F";width 2.693in;height 2.0781in;depth
0in;original-width 6.4999in;original-height 5.0004in;cropleft "0";croptop
"1";cropright "1";cropbottom "0";filename 'Fig3.jpg';file-properties
"XNPEU";}} \\ 
Fig. 3: solutions for $\beta =10000,N=40$%
\end{tabular}%
\end{array}%
\end{equation*}

For $N=64,$ $\Delta t=0.000005,$ $\beta =10000,$ on the interval $%
[-0.2,1.06] $ results aredocumented in Table 2 to compare with Ref \cite{fs5}%
, in terms of $L_{\infty }$ norm, at different time steps.%
\begin{equation*}
\begin{tabular}{lclll}
\multicolumn{5}{l}{\textbf{Table 2}} \\ 
\multicolumn{5}{l}{$L_{\infty }$ errors norms at some different times for
the first test problem for $N=64$} \\ \hline
Method & $t=0.0005$ & $t=0.0015$ & $t=0.0025$ & $t=0.0035$ \\ \hline
Present & $1.02\times 10^{-2}$ & $1.49\times 10^{-1}$ & $3.24\times 10^{-1}$
& $4.78\times 10^{-1}$ \\ 
\cite{fs5} & $2.55\times 10^{-3}$ & $1.62\times 10^{-2}$ & $8.65\times
10^{-2}$ & $6.98\times 10^{-2}$ \\ \hline
\end{tabular}%
\end{equation*}

\textbf{(b) }Secondly, the initial pulse profile%
\begin{equation*}
u(x,0)=\sec h^{2}(10x)
\end{equation*}%
is chosen as the initial condition for our first numerical experiment. Under
the boundary conditions (\ref{s3}) we obtained the solutions. In numerical
calculations, the constants in Eq.(\ref{f2}) are selected $\lambda =0.1$ and 
$\beta =1.$ For the discretization of space and time, space/time incerements 
$h=0.005$ and $\Delta t=0.05$ are used. Then, the algorithm is run up to
time $t=40$ over the domain $[-50,50].$ In short period Figure 4 is drawn
for $[a,b]=[-2,2]$ at times $t=0.1,$ $0.2,$ $0.3,$ $0.4$ and $0.5$ in which
we see that the diffusion is more dominat than the reaction. Figure 5 is
drawn for $[a,b]=[-6,6]$ at $t=0,$ $1,$ $2,$ $3,$ $4,$and $5$. \ After the
concentration reached the lowest level, it start to increase up to a level $%
u=1,$so that the reaction are observed over the the diffusion. Finally, in
Figure 6, solutions are depicted at $t=5,$ $10,15,$ $20,$ $25$, $30,$ $35$
and $40$ in the space interval $[a,b]=[-50,50].$As time advance more, the
concentration return to initial form, then flatten through both sides with
sharp lateral fronts. Thus solutions looks like bell-shape with flat top. We
see \ that the diffusion is tottaly efficient in advance time. 
\begin{equation*}
\begin{array}{c}
\begin{array}{cc}
\begin{tabular}{l}
\FRAME{itbpF}{2.693in}{2.111in}{0in}{}{}{fig4.jpg}{\special{language
"Scientific Word";type "GRAPHIC";maintain-aspect-ratio TRUE;display
"USEDEF";valid_file "F";width 2.693in;height 2.111in;depth
0in;original-width 6.4999in;original-height 5.0834in;cropleft "0";croptop
"1";cropright "1";cropbottom "0";filename 'Fig4.jpg';file-properties
"XNPEU";}} \\ 
Fig. 4: Solutions at early times. \\ 
with at $t=0.1,$ $0.2,$ $0.3,$ $0.4,$ $0.5.$%
\end{tabular}
& 
\begin{tabular}{l}
\FRAME{itbpF}{2.693in}{2.0954in}{0in}{}{}{fig5.jpg}{\special{language
"Scientific Word";type "GRAPHIC";maintain-aspect-ratio TRUE;display
"USEDEF";valid_file "F";width 2.693in;height 2.0954in;depth
0in;original-width 6.4999in;original-height 5.0436in;cropleft "0";croptop
"1";cropright "1";cropbottom "0";filename 'Fig5.jpg';file-properties
"XNPEU";}} \\ 
Fig. 5: Short-time behavior. \\ 
with at $t=0,1,$ $2,$ $3,$ $4,$ $5.$%
\end{tabular}%
\end{array}
\\ 
\begin{tabular}{l}
\FRAME{itbpF}{2.693in}{2.0773in}{0in}{}{}{fig6.jpg}{\special{language
"Scientific Word";type "GRAPHIC";maintain-aspect-ratio TRUE;display
"USEDEF";valid_file "F";width 2.693in;height 2.0773in;depth
0in;original-width 6.4999in;original-height 5.0004in;cropleft "0";croptop
"1";cropright "1";cropbottom "0";filename 'Fig6.jpg';file-properties
"XNPEU";}} \\ 
Fig. 6: Long-time behavior. \\ 
with at $t=0,$ $5,$ $10,$ $15,$ $20,$ $25,$ $30,$ $35,.40.$%
\end{tabular}%
\end{array}%
\end{equation*}%
In the Table 3 errors at some different times for the second test problem is
shown.

\begin{equation*}
\begin{tabular}{lcllll}
\multicolumn{6}{l}{\textbf{Table 3:}$\text{Relative e}$rrors at some
different times for the second test problem b for $N=64$} \\ 
\multicolumn{6}{l}{Parameters: $\lambda =0.1,$ $\beta =1,$ $\Delta t=0.05$
and $x\in \lbrack -50,50]$} \\ \hline
Method & $t=5$ & $t=10$ & $t=15$ & $t=20$ & $t=40$ \\ \hline
Present & $1.383\times 10^{-2}$ & $7.834\times 10^{-3}$ & $6.029\times
10^{-3}$ & $5.066\times 10^{-3}$ & $3.416\times 10^{-3}$ \\ 
\cite{fs5} & $1.386\times 10^{-2}$ & $7.860\times 10^{-3}$ & $6.054\times
10^{-3}$ & $5.090\times 10^{-3}$ & $3.434\times 10^{-3}$ \\ \hline
\end{tabular}%
\end{equation*}

\section{Conclusion}

The collocation methods with trigonometric B-spline functions is made up to
find solutions the Fisher's equation. We have shown that method is capable
of producing solutions of the Fisher's equation fairly. The method can be
used as an alternative to the more usual assosiate B-spline collocation and
Galerkin methods.

\textbf{Acknowledgement:} This paper have been presented at the
International Conference on Natural Science and Engineering, 2016, Kilis,
Turkey.


\begin{thebibliography}{99}
\bibitem{first} R.A. Fisher, "The Wave of Advance of Advantageous Genes", 
\textit{Annals of Human Genetics,} vol. 7, no.4, pp. 355-369, 1973.

\bibitem{conv} T. Mavoungou and Y. Cherruault, "Numerical Study of Fisher's
Equation by Adomian's Method", \textit{Mathematical and Computer Modelling},
vol. 19, no. 1, pp. 89-95, 1994.

\bibitem{g1} A. Nikolis, "Numerical Solutions of Ordinary Differential
Equations with Quadratic Trigonometric Splines", \textit{Applied Mathematics
E-Notes}, vol. 4, pp. 142-149, 1995.

\bibitem{g2} A. Nikolis and I. Seimenis, "Solving Dynamical Systems with
Cubic Trigonometric Splines", \textit{Applied Mathematics E-notes,} vol. 5,
pp. 116-123, 2005.

\bibitem{we} Nur Nadiah Abd Hamid , Ahmad Abd. Majid, and Ahmad Izani Md.
Ismail, "Cubic Trigonometric B-Spline Applied to Linear Two-Point Boundary
Value Problems of Order Two", \textit{World Academy of Science, Engineering
and Technology}, vol. 4, no. 10, pp. 1377-1382, 2010.

\bibitem{w1} Y. Gupta \ and M. Kumar, "A Computer based Numerical Method for
Singular Boundary Value Problems", \textit{International Journal of Computer
Applications,} vol. 30, no. 1, pp. 21-25, 2011.

\bibitem{aa} Abbas, A. A. Majid, A. I. Md. \.{I}smail and A. Rashid, "The
Application of Cubic Trigonometric B-spline to the Numerical Solution of the
Hyperbolic Problems", \textit{Applied Mathematica and Computation}, vol.
239, pp. 74-88, 2014.

\bibitem{ma} Abbas, A. A. Majid, A. I. Md. \.{I}smail and A. Rashid,
"Numerical Method Using Cubic Trigonometric B-spline Technique for
non-classical Diffusion Problems", \textit{Abstract and applied analysis, }%
2014

\bibitem{rubin} S. G. Rubin and R. A. Graves, Cubic Spline Approximation for
Problems in Fluid Mechanics, Nasa TR R-436, Washington, DC, 1975.

\bibitem{k1} S. Li, L. Petzold and Y. Ren, (1998). "Stability of Moving Mesh
Systems of Partial Differential Equations", \textit{SIAM Journal on
Scientific Computing}, 20, 719--738.

\bibitem{k2} Y. Qiu and D. M. Sloan, "Numerical solution of Fisher's
equation using a moving mesh method", \textit{J. Comput. Phys.}, vol. 146,
pp. 726--746, 1998.

\bibitem{z1} M. Ablowitz and A. Zepetella, "Explicit Solution of Fisher's
Equation for a Special Wave Speed", \textit{Bulletin of Mathematical Biology}%
, vol. 41, pp. 835--840, 1979.

\bibitem{fs9} Danumjaya P and Pani AK. "Orthogonal cubic spline collocation
method for the extended Fisher-Kolmogorov equation". \textit{Journal of
Comput. and Appl. Math}., vol. 174, pp. 101-117, 2015.

\bibitem{fs4} A. Sahin, I. Dag and B. Saka, "A B-spline algorithm for the
numerical solution of Fisher's equation", \textit{Kybernetes,} vol. 37, pp.
326--342, 2008.

\bibitem{fs5} \.{I} Dag, A. \c{S}ahin and A. Korkmaz, "Numerical
Investigation of the Solution of Fisher's Equation via the B-Spline Galerkin
Method", \textit{Numer Methods Partial Differential Eq.} vol. 26, pp.
1483--1503, 2010.

\bibitem{fs8} R. C. Mittal and G. Arora, "Quintic B-spline collocation
method for numerical solution of the Extended Fisher -Kolmogorov equation", 
\textit{Int. J. of Appl. Math and Mech}., vol. 6, no. 1, pp. 74-85, 2010.

\bibitem{fs3} RC. Mittal, G. Arora, "Efficient numerical solution of
Fisher's equation by using B-spline method". \textit{Int. J. Comput. Math}.
vol. 87, no. 13, pp. 3039--3051, 2010.

\bibitem{fs6} C. G. Zhu, , W. S. Kang, "Numerical solution of
Burgers--Fisher equation by cubic B-spline quasi-interpolation", \textit{%
Applied Mathematics and Computation}, vol. 216, no. 9, pp. 2679--2686, 2010.

\bibitem{fs2} R. C. Mittal and R. K. Jain, "Numerical solutions of nonlinear
Fisher's reaction--diffusion equation with modified cubic B-spline
collocation method", \textit{Mathematical Sciences}, vol. 7, no. 12, 2013.

\bibitem{fs7} Ali Sahin and \"{O}zlem \  \"{O}zmen, "Usage of Higher Order
B-splines in Numerical Solution of Fisher's Equation", \textit{International
Journal of Nonlinear Science}, vol. 17, no. 3, pp. 241-253, 2014.

\bibitem{fs1} M. Aghamohamadi, J. Rashidinia, R. Ezzati, "Tension spline
method for solution of non-linear Fisher equation", \textit{Applied
Mathematics and Computation} vol. 249, pp. 399--407, 2014.

\bibitem{fs0} O. Ersoy and I. Dag, "The Extended B-spline Collocation Method
for Numerical Solutions of Fisher's equation", \textit{AIP Conference
Proceedings} vol. 1648, no. 370011, 2015.

\bibitem{fs10} I Dag and O Ersoy, "The exponential Cubic B-spline Algorithm
for Fisher Equation", C\textit{haos, Solitons and Fractals} vol. 86, pp.
101--106, 2016.
\end{thebibliography}
\end{document}